\documentclass[12pt,a4paper]{article}
\pdfoutput=1
\usepackage{setspace}
\usepackage[english]{babel}		% Para a hifenização em português
\usepackage[latin1]{inputenc}
\usepackage{amscd}
\usepackage{amsfonts}
\usepackage{amsgen}
\usepackage{amsmath}
\usepackage{amssymb}
\usepackage{amstext}
\usepackage{amsthm}
\usepackage{graphicx}
\usepackage{latexsym}
\usepackage{makeidx}																% Para usar Índice Remissivo
\usepackage{epsfig}
\usepackage{wrapfig}
\usepackage[all,knot,arc,import,poly]{xy}
\usepackage[usenames]{color}
\usepackage[mathscr]{eucal}
\usepackage{indentfirst}
\usepackage[T1]{fontenc}
\usepackage{setspace}
\newtheorem{Def}{Definition}[section]
\newtheorem{Teo}[Def]{Theorem}
\newtheorem{Prop}[Def]{Proposition}
\newtheorem{Lema}[Def]{Lemma}

\newcommand{\MS}{\mathcal{S}}

\newcommand{\MC}{\mathcal{C}}
\newcommand{\menos}{\backslash}

\newcommand{\setad}{\rightarrow}

\newcommand{\sse}{\Leftrightarrow}

\newcommand{\sigx}{\sigma(x)}
\newcommand{\sigy}{\sigma(y)}
\newcommand{\sigxy}{\sigma(xy)}

\newcommand{\rgs}{(RG)^+}

\newcommand{\uni}{\mathcal U}
\newcommand{\centro}{\mathcal Z}

\newcommand{\inv}{^{-1}}

\newcommand{\dem}{\begin{proof}}
\newcommand{\cqd}{\end{proof}}

\newcommand{\che}{\left\{}
\newcommand{\chd}{\right\}}
\newcommand{\pare}{\left(}
\newcommand{\pard}{\right)}
\newcommand{\lan}{\langle}
\newcommand{\ran}{\rangle}

	\newcommand{\properpagestyle}
	
\makeindex

\title{Anticommutativity of Symmetric Elements under Generalized Oriented Involutions}
\author{Tonucci, L. E. \and Petit Lobão, T. C.}

\begin{document}
\maketitle
%\input{enrolacoes}

%%%%%%%%%%% CAPA %%%%%%%%%%
%\include{Capa}

%----------------------------------INDICE------------------------------
\renewcommand{\contentsname}{Índice}
%\tableofcontents \mainmatter
%---------------------------------------------------------------------

%\newpage
%\pagestyle{myheadings}
\pagenumbering{arabic}
\setcounter{page}{1}
%
%\def\nonumchapter#1{%
%    \chapter*{#1}
%    \addcontentsline{toc}{chapter}{#1}}
%
%{\properpagestyle \begin{titlepage}\include{Introducao}\end{titlepage}}%introdução

%\addcontentsline{toc}{chapter}{Introdução}

\section{Abstract}

Let $R$ be a ring with $char(R)\neq2$ whose unit group are denoted by $\uni(R)$, $G$ a group with involution $*$, and $\sigma:G\setad\uni(R)$ a nontrivial group homomorphism, with $ker\ \sigma=N$, satisfying $xx^*\in N$ for all $x\in G$. Let $RG$ be the group ring of $G$ over $R$ and define the involution $\sigma*$ in $RG$ by $\pare \sum_{x\in G}\alpha_xx\pard^{\sigma*}=\sum_{x\in G}\sigx\alpha_xx^*$. In this paper, we will classify the group rings $RG$ such that $\MS$ is anticommutative, where $\MS$ is the largest subset of $(RG)^+=\che \alpha\in RG: \alpha^{\sigma *}=\alpha\chd$ that can satisfy anticommutativity under $char(R)\neq2$.
\vspace{.1cm}\\
{\it Keywords:} group rings, rings with involution, symmetric elements, generalized oriented involution.

\section{Introduction}

Let $RG$ be a group ring of a group $G$ over a commutative ring $R$ with unity. Given $*$ an involution in $G$, we can naturally induce an involution in $RG$, defined by the linear extension of $*$. 

In the context of $K$-theory, Novikov in \cite{N70} introduced a really interesting nonlinear involution, using a nontrivial homomorphism $\sigma:G\setad\che\pm1\chd$ such that $\sigma(x^*x\inv)\in ker\sigma$, denoted by $\sigma*$, which maps $\sum_{x\in G}\alpha_xx\mapsto \sum_{x\in G}\alpha_x\sigx x^*$, the so called oriented involution.

Given a ring $R$ with involution $*$, we can define the set of symmetric elements, in respect of this involution, given by $R^+=\che\alpha\in R:\alpha^*=\alpha\chd.$
In a fundamental work, \cite{A68}, Amistur proved that, if $R^+$ satisfies a polynomial identity, then so does $R$. Thus, asking about which identities in $R^+$ could be lifted to $R$ is quite natural, and since then, many authors turn their attention to this question.

Interesting identities can be given by the Lie bracket, $[\alpha,\beta]=\alpha\beta-\beta\alpha$, and the Jordan operator, $\alpha\circ\beta=\alpha\beta+\beta\alpha$. We say that $R^+$ is commutative, if the restriction of the Lie bracket to it is identically $0$, in the same way, we say that it is anticommutative if the restriction of the Jordan product is $0$. In \cite{JM06,BP06,GP13a,GP14}, the authors classify the group rings such that $\rgs$ is commutative or anticommutative. Generalizing these identities, the papers \cite{CP12,LSS09,GPS09,L99,L00} classify $RG$ when $\rgs$ satisfies Lie nilpotency, in other words, when $[x_1,\ldots,x_n]=0,\ \forall x_i\in \rgs$, and Lie $n$-Engel, when $[x,y,y,\ldots,y]=0,\ \forall x,y\in \rgs$.

Aiming for generalization, replacing the orientation $\sigma:G\setad\che\pm1\chd$ with an homomorphism, which we denote similarly, $\sigma:G\setad\uni(R)$, we can also induce an involution in $RG$, if the involution $*$ is compatible with $\sigma$ in the sense that $xx^*\in N=ker\sigma,\ \forall x\in G$. This new involution is becoming a researching object, and some of the identities explicited above were studied in \cite{V13}.

In this paper, we will classify the group rings $RG$, such that a suitable subset of $\rgs$ is anticommutative, in other words, in which the Jordan product is trivial, where $\rgs$ are the set of symmetric elements under a generalized oriented involution. Since in $char(R)=2$, the Lie bracket coincides with the Jordan product, we will assume that $char(R)\neq2$ and will use this fact without further mention.

The kernel of $\sigma$ is denoted by $N$, the symmetric elements under $*$ in $G$ are collected in $G_*$, the center of $G$, in $\centro(G)$, $(x,y)=x\inv y\inv xy$ is the multiplicative commutator, and $x^y=y\inv xy$ is the conjugation of $x$ by $y$.

\section{Symmetric Elements Anticommute}

Firstly, we will obtain a set of generators of the symmetric elements $(RG)^+$.

Let $N_*=G_*\cap N$. For $\alpha\in RG$, write
$$\alpha=\sum_{x\in N_*}\alpha_xx+\sum_{x\in N\menos G_*}\alpha_xx+\sum_{x\in G_*\menos N}\alpha_xx+\sum_{x\in G\menos(G_*\cup N)}\alpha_xx,$$
then
$$\alpha^{\sigma*}=\sum_{x\in N_*}\alpha_x x+\sum_{x\in N\menos G_*}\alpha_x x^*+\sum_{x\in G_*\menos N}\sigx\alpha_x x+\sum_{x\in G\menos(G_*\cup N)}\sigx\alpha_x x^*,$$
so
$$\alpha^{\sigma*}=\alpha\ \text{if and only if}\che
\begin{array}{ll}
\alpha_x=\alpha_{x^*}&\text{for}\ g\in N\menos G_*\\
\alpha_x=\sigx\alpha_{x}&\text{for}\ g\in G_*\menos N\\
\alpha_x=\sigx\alpha_{x^*}&\text{for}\ g\in G\menos (G_*\cup N),
\end{array}\right.$$
thus $(RG)^+$ is spanned over $R$ by elements in the sets,
$$\begin{array}{l}
\MS_1=\che x:x\in N_*\chd,\\
\MS_2=\che \alpha x:x\in G_*\menos N,\ \alpha(1-\sigx)=0\chd,\\
\MS_3=\che x+\sigx x^*: x\in G\menos G_*\chd.
%\MS_3=\che x+\sigx x^*: x\in N\menos G_*\chd,\\
%\MS_4=\che x+\sigx x^*: x\notin (G_*\cup N)\chd.
\end{array}$$

The elements in $(RG)^+$ anticommute if and only if any two in the union of the $\MS_i$ anticommute. Note that, if $x\in N_*$, then $x\in \MS_1$; thus if $(RG)^+$ is anticommutative, then $2x^2=x^2+x^2=0$, and $char(R)=2$. Since we are avoiding $char(R)=2$, in fact, we will not study when $(RG)^+$ is anticommutative, but the largest subset of $(RG)^+$ that can satisfy this idendity under $char(R)\neq2$. In order to do this, we will replace $\MS_1$ with $2\MS_1=\che 2x:x\in N_*\chd$ and study the set $\MS$ given by the spanning over the set $2\MS_1\cup \MS_2\cup\MS_3$.

%Given $x\in G_*$, due the compatibility condition, $\sigx^2=\sigma(x^*x)=1,$ so $x^2\in N$, then taking $\alpha=1+\sigx$, we conclude $\alpha(1-\sigx)=0$, then $\MS$ could be rewritten by the spanning of
%$$\che \alpha x:x\in G_*\menos N,\ \alpha(1-\sigx)=0\chd\cup\che x+\sigx x^*: x\in G\chd.$$

In the case $\sigma=\che\pm1\chd$, we have the following result, proved in \cite{GP14}, which will be very helpfull.

\begin{Teo}[Theorem 2.2, \cite{GP14}]Let $x\mapsto x^*$ denote an involution on a group $G$, and let $\sigma:G\setad\che\pm1\chd$ be a nontrivial orientation homomorphism and compatible with $*$ in the sense that $\sigma(x^*)=\sigma(x)$ for all $x\in G$. Let $R$ be a commutative ring with unity and of characteristic different from $2$. Let $R_2$ be the set of $r\in R$ satisfying $2r=0$. For $\alpha=\sum_{x\in G}\alpha_xx$ in the group ring $RG$, define $\alpha^{\sigma*}=\sum_{x\in G}\sigx\alpha_x x^*$. If the elements of
$$\mathcal{C}=\che x+\sigx x^*: x\in N_*\cup (G\menos G_*)\chd\cup\che\alpha x:x\in G_*\menos N, \alpha\in R_2\chd$$
anticommute, then either one of the following statements holds:
\begin{itemize}
\item [(A)] $char(R)=4$ or $8$, $G$ is abelian and $*=Id$, or
\item [(B)] $char(R)=4$, $G$ is abelian and $*=Id|_N$, or
\item [(C)] $char(R)=4$,
\begin{itemize}
\item [(i)] $G'=\che 1,s\chd$ for some $s\neq 1$ (implicity $s^2=1$),
\item [(ii)] $x^*\in \che x,sx\chd$, for any $x\in G$,
\item [(iii)] either the $*$-symmetric elements not in $N$ commute or $R_2^2=\che 0\chd$.
\end{itemize}
\end{itemize}
Conversely, if $G$ is a group with an index two subgroup $N$ and $\sigma:G\setad\che\pm1\chd$ is the orientation homomorphism with kernel $N$, then the elements of $\mathcal{C}$ anticommute in any of the three specified situations.
\label{teo.gp14}
\end{Teo}

Note that $C=\che x\in G:\sigx=\pm1\chd$, due to the compatibility of $*$ and $\sigma$, is a $*$-invariant subgroup of $G$. So, as $\MC$ is an restriction of $\MS$, if $\MS$ is anticommutative, then so does $\MC$, thus $C\leq G$ satisfies one of the above conditions.

\begin{Lema} If $\MS$ is anticommutative, then $char(R)=4$ or $8$.\label{as.char4.char8}\end{Lema}

\dem Notice that $1\in N_*$, hence $2$ anticommutes with itself, so, $4=2^2=-2^2=-4$, thus, $8=0$, which implies $char(R)=4$ or $8$, since $char(R)\neq2$.\cqd

%\begin{Prop} Se $*=Id$, então $\MS$ é anticomutativo se, e somente se, $2(1+\sigx+\sigy+\sigxy)=0,\ \forall x,y\in G$.
%\label{as.*=id}%
%\end{Prop}

%\dem Observe que se $*=Id$, então $xy=x^*y^*=(yx)^*=yx,\ \forall x,y\in G$, logo $G$ é abeliano, portanto 
%$$\begin{array}{c}
 %\MS\ \text{é anticomutativo}\\
 %\Updownarrow\\
  % (x+\sigx x^*)(y+\sigy y^*)+(y+\sigy y^*)(x+\sigx x^*)=0,\ \forall x,y\in G\\
   %\Updownarrow\\
 %2(1+\sigx+\sigy+\sigxy)xy=0\,\ \forall x,y\in G\\
 %\Updownarrow\\
 %2(1+\sigx+\sigy+\sigxy)=0,\ \forall x,y\in G.
 %\end{array}$$\cqd

\begin{Lema}Suppose that $\MS$ is anticommutative. If $*\neq Id$, then $char(R)=4$, $xx^*=x^*x$, and $x^2\in G_*$, for all $x\notin G_*$.
\label{as.*.neq.id}
\end{Lema}

\dem  Let $x\notin G_*$, by hypothesis, $(x+\sigx x^*)^2+(x+\sigx x^*)^2=0$, so
$$2(x^2+\sigx xx^*+\sigx x^*x+\sigx^2(x^*)^2)=0.$$
Since the left member of the equation is $0$, and $2\neq0$, we must have $x^2\in\che xx^*,x^*x,(x^*)^2\chd$. As $x\neq x^*$, the only possibility is $x^2=(x^*)^2=(x^2)^*$, so, $x^2\in G_*$. For the same reason, we can find $xx^*=x^*x$, and, for the equation holds, it is necessary that $4\sigx xx^*=0$, witch implies $4=0$, since $\sigx\in\uni(R)$ and $RG$ is freely generated by $G$.\cqd

From now on, the above lemma will be used without any reference.

Assuming that $x,y\in G$ and $\MS$ is anticommutative, we have the following equation that will be used several times in this work:
\begin{equation}
 \begin{array}{rcl}
0&=&(x+\sigma(x)x^*)(y+\sigma(y)y^*)+(y+\sigma(y)y^*)(x+\sigma(x)x^*)\\
&=&xy+yx+\sigma(xy)x^*y^*+\sigma(xy)y^*x^*+\sigma(y)xy^*+\sigma(y)y^*x+\sigma(x)yx^*+\sigma(x)x^*y.
\end{array}
\label{as.eq.prin}
\end{equation}

\begin{Lema} Suppose that $\MS$ is anticommutative. Given $x,y\in G$, then, $xy=yx$ if, and only if, $x^*y=yx^*$. Moreover, if $x,y\notin G_*$, then $xy=yx=x^*y^*=y^*x^*$, $xy^*=y^*x=x^*y=yx^*$, and $2(1+\sigxy)=2(\sigx+\sigy)=0$.\label{as.comut}\end{Lema}

\dem  If $x\in G_*$, trivialy we have the equivalence. If $y\in G_*$, then, apllying the involution in both sides of $xy=yx$, we get that $yx^*=x^*y$, and the result holds.

We can assume that $x,y\notin G_*$. Applying the involution to $xy=yx$, we get $y^*x^*=x^*y^*$, so equation (\ref{as.eq.prin}) could be written as
$$2xy+2\sigma(xy)x^*y^*+\sigma(y)xy^*+\sigma(y)y^*x+\sigma(x)yx^*+\sigma(x)x^*y=0.$$

As $char(R)\neq2$, we must have $xy=yx\in\che x^*y^*,xy^*,y^*x,yx^*,x^*y\chd$; since $x,y\notin G_*$, $xy=yx\notin\che xy^*,y^*x,yx^*,x^*y\chd$, so $xy=yx=x^*y^*=y^*x^*$. Replacing this in the equation above, we get $2(1+\sigxy)=0$. Hence
\begin{equation}                                                                                                                                                                                                   
\sigma(y)xy^*+\sigma(y)y^*x+\sigma(x)yx^*+\sigma(x)x^*y=0,                                                                                                                                                                                                      
\label{as.eq.1}
\end{equation}
and $x^*y$ must be equal to another of the three remaining elements in the support of the left side of this equation. 

If $x^*y=y^*x$, (\ref{as.eq.1}) gives $xy^*=yx^*$, thus, $x(y^*y)=yx^*y=yy^*x=(y^*y)x$, and, as $(x,y)=1=(x,y^*y)$, hence $(x,y^*)=1$; if $x^*y=xy^*$, in the same way we find $(x,y^*)=1$. Replacing $yx^*$ by $x^*y$       in (\ref{as.eq.1}), since $2\sigx\neq0$, we must have $xy^*=y^*x=x^*y=yx^*$, and $2(\sigx+\sigy)=0$.

To the converse, it is sufficient to apply the result to $x^*$ and $y$.
\cqd

\begin{Lema}Suppose that $\MS$ is anticommutative. If $x,y \notin G_*$ and $\sigy\neq-1$, then, $x^y\in\che x^*,x\chd$ or, $\sigxy=-1$, $xy=x^*y^*$ and $xy^*=x^*y$.
\label{as.base.1}\end{Lema}

\dem Suppose that $x^y\notin \che x,x^*\chd$. Since the equation (\ref{as.eq.prin}) holds, and $x,y\notin G_*$ we have that $xy\in\che x^*y^*,y^*x^*,y^*x\chd$. Let us study these three possibilities.

If $xy=y^*x$, as $\sigy\neq-1$, (\ref{as.eq.prin}) holds only if $xy$ is equal to another element in its support. As the unique possibilities are $xy=x^*y^*$ and $xy=y^*x^*$, we have that $xy=y^*x=x^*y^*$, since $x\notin G_*$, it implies $xy=y^*x\neq y^*x^*$. As $xy=y^*x=x^*y^*\notin \che yx,y^*x^*,x^*y,xy^*,yx^*\chd$, then, applying the involution, $y^*x^*=x^*y=yx\notin \che x^*y^*,xy,y^*x,yx^*,xy^*\chd$, and replacing it in (\ref{as.eq.prin}), we get
$$(1+\sigy+\sigxy)xy+(1+\sigx+\sigxy)yx+\sigx yx^*+\sigy xy^*=0,$$
thus $yx^*=xy^*$ and
$$(1+\sigy+\sigxy)=(1+\sigx+\sigxy)=(\sigx+\sigy)=0,$$
which gives a contradiction, since $2\sigx\neq0$.

If $xy=y^*x^*$, then $xy\neq x^*y^*$, since $y^*x^*=x^*y^*$, applying the involution, it implies $xy=yx$, a contradiction. By (\ref{as.eq.prin}), we have $\sigxy=-1$, then $2xy\in\MS$, so it anticommutes with $(x+\sigx x^*)$, in the other words,
 $$2xyx+2\sigx xyx^*+2x^2y+2\sigx xx^*y=0,$$
and, since $(x,y)\neq1$ and $x\notin G_*$, $x^2y=xyx^*$; which implies $xy=yx^*$, a contradiction.

If $xy=x^*y^*$ and $\sigxy\neq-1$, then $xy\in\che yx,y^*x^*,x^*y,yx^*,xy^*,y^*x\chd$, which leads to a contradiction by hypothesis or to the cases above. So $\sigxy=-1$. Applying the involution in $xy=x^*y^*$, we get $yx=y^*x^*$, and (\ref{as.eq.prin}) implies $xy^*\in\che x^*y,yx^*,y^*x\chd$. As $xy\neq yx$, then $xy^*\neq y^*x$, so $xy^*=yx^*$ or $xy^*=x^*y$. Suppose that $xy^*=yx^*$, thus, $x^*y=y^*x$. As $\sigxy=-1$, then $\sigx=-\sigy\inv$; on the other hand, $\sigy\neq-1$, so, $x\notin N$, hence $\sigma(x^2y)\neq-1$.
 
 Notice that $x^2y\notin G_*$, since $(x^2y)^*=y^*x^*x^*=yxx^*=yx^*x=xy^*x=xx^*y\neq x^2y$; then, applying what we proved above to $x$ and $x^2y$, we obtain that $x^{x^2y}\in \che x,x^*\chd$ or $x(x^2y)=x^*(x^2y)^*$. As $x^{x^2y}=x^y$, we can assume that $x(x^2y)=x^*(x^2y)^*$, otherwise $xy\in\che yx,y^*x^*,x^*y,yx^*,xy^*,y^*x\chd$, which is false. So $x^3y=x^*y^*(x^2)^*=xy(x^2)^*=xyx^2$, thus $x^2y=yx^2$. On the other hand, $x^3y=x^*y^*x^*x^*=x^*yxx^*=y^*xxx^*=x^2y^*x^*$, implying $xy=y^*x^*=yx$, a contradiction. Then $xy^*=x^*y$.
\cqd

\begin{Lema}Suppose $\MS$ is anticommutative. If $x,y \notin G_*$ and $\sigy\neq-1$, then only one of the following holds:
\begin{itemize}
 \item [(i)] $xy=yx=x^*y^*=y^*x^*$ and $2(1+\sigxy)=0=2(\sigx+\sigy)$.
 \item [(ii)] $xy=yx^*=y^*x=x^*y^*$ and $1+\sigx+\sigy+\sigxy=0$.
 \item [(iii)] $xy=x^*y^*\neq yx^*=y^*x$, $\sigxy=-1$, and $\sigx=-\sigy$.
 \item [(iv)] $xy=yx^*\neq x^*y^*=y^*x$ and $\sigx=-1$.
\end{itemize}
\label{as.base.2}\end{Lema}

\dem If $(x,y)=1$, Lemma \ref{as.comut} gives (i). 

Suppose $(x,y)\neq1$, then, by Lemma \ref{as.base.1}, $xy=yx^*$, or $xy=x^*y^*$, $x^*y=xy^*$, and $\sigxy=-1$. If $xy=yx^*=x^*y^*$, using a similar argument to the case $xy=y^*x$ of the lemma above, we can find $xy=yx^*=x^*y^*=y^*x$; so, applying the involution and replacing it in equation (\ref{as.eq.prin}), we get (ii).

If $xy\neq yx^*$, then $xy=x^*y^*$, $\sigxy=-1$ and $xy^*=x^*y$, thus $y^*x^*=yx$ and $yx^*=y^*x$. Since $(x,y)\neq1$ and $x\notin G_*$, $yx\neq yx^*\neq x^*y$, so $yx^*\notin\che xy,yx,x^*y^*,y^*x^*,xy^*,x^*y\chd$; then (\ref{as.eq.prin}) implies $\sigx=-\sigy$ and (iii) holds.

Finally, suppose that $xy=yx^*\neq x^*y^*$. Untill now, we have already proved that, under the lemma assumptions, it holds (i)-(iii) or $xy=yx^*\neq x^*y^*$; so appying it to $x^*$ and $y$, we get one of the above:
\begin{itemize}
 \item [(a)] $x^*y=yx^*=xy^*=y^*x$.
 \item [(b)]  $x^*y=yx=y^*x^*=xy^*$
 \item [(c)] $x^*y=xy^*\neq yx=y^*x^*$.
 \item [(d)] $x^*y=yx\neq xy^*$.
\end{itemize}
Notice that (a)-(c) are equivalent to (i)-(iii), so they not occur, then $x^*y=yx$; so, applying the involution, $y^*x=x^*y^*$ and replacing it in (\ref{as.eq.prin}), we find $\sigx=-1$.
\cqd

\begin{Lema}Suppose that $\MS$ is anticommutative. If $x,y\notin G_*$, then the following holds:
\begin{itemize}
 \item [(i)] $xy\in\che yx, yx^*,y^*x,x^*y^*\chd$.
 \item [(ii)] $xy=yx$ if, and only if, $xy\in G_*$.
 \item [(iii)] $xy=yx^*$ if, and only if, $x^*y=yx$.
 \end{itemize}
\label{as.2.n.simet}\end{Lema}

\dem To prove (i), suppose $(x,y)\neq1$. If $\sigy\neq-1$, then, by Lemma \ref{as.base.1}, we have the result. If $\sigx\neq-1$, applying Lemma \ref{as.base.1} to $y^*$ and $x$, we get $y^*x=xy$ or $(y^*)x^*=(y^*)^*x=yx$, and both cases implies (i).

Then we can suppose $\sigx=\sigy=-1$, hence $(x-x^*)$ anticommutes with $(y-y^*)$, thus
$$xy+yx+x^*y^*+y^*x^*-x^*y-yx^*-y^*x-xy^*=0,$$
so, $xy\in\che yx,x^*y^*,y^*x^*,x^*y,yx^*,y^*x,xy^*\chd$.
By hyphotesis, $xy\notin\che yx,x^*y,xy^*\chd$. If $xy\notin\che x^*y^*,yx^*,y^*x\chd$, then $xy=y^*x^*$; but, as $char(R)=4$, $2xy$ can not be canceled with another element in equation (\ref{as.eq.prin}); which leads to a contradiction, so, $xy\in\che x^*y^*,yx^*,y^*x\chd$ and (i) holds.

To prove (ii), suppose $xy\in G_*$, so $xy=y^*x^*$. By (i), $xy\in\che yx^*,y^*x,x^*y^*\chd$, but as $xy=y^*x^*$ and $x,y\notin G_*$, it follows that $xy\notin\che yx^*,y^*x\chd$, so, $y^*x^*=xy=x^*y^*$, which implies $(x,y)=1$. The converse is consequence of Lemma \ref{as.comut}.

To verify (iii), suppose $xy=yx^*$ and $x^*y\neq yx$. In this case, $(x,y)\neq1\neq (x^*,y)$, so, applying item (i) to $x^*$ and $y$ we find $x^*y\in\che y^*x^*,xy^*\chd$. If $x^*y=y^*x^*$, then $y^*x=xy=yx^*$, so, $x^*y=y^*x^*=xy^*$; if $x^*y=xy^*$, then $y^*x=yx^*=xy$, so, $x^*y=xy^*=y^*x^*$. Thus, in both cases, $x^*y=y^*x^*=xy^*\neq yx$, and applying the involution, we get $y^*x=xy=yx^*\neq x^*y^*$, so, by (\ref{as.eq.prin}), 
$$\sigxy+\sigx+\sigy=1+\sigx+\sigy=1+\sigxy=0,$$
which leads to a contradiction, since $char(R)\neq2$. To the converse, it is enough to apply the same result to $x^*$ and $y$.
\cqd

\begin{Lema}If $\MS$ is anticommutative, so $(x^2,y)=1,\ \forall x,y\notin G_*$.\label{as.quad.n.simet.comut}\end{Lema}

\dem Let $x,y\notin G_*$ such that $(x,y)\neq1$, consequently, by item (ii) in Lemma \ref{as.2.n.simet}, $xy\notin G_*$. If $x,y\in C$, by Theorem \ref{teo.gp14}, the result holds.

If $x\in C$ and $y\notin C$, then, applying Lemma \ref{as.base.2} to $x$ and $y$, we get $xy=yx^*$, then $x^2y=xyx^*=y(x^*)^2=yx^2$. If $x\notin C$ and $y\in C$, we can procede in an analogous way.

If $x,y\notin C$, by the same lemma, (ii) or (iii) occurs. If (ii) holds, the result follows as above. If (iii) holds, then applying Lemma \ref{as.base.2} to $xy$ and $x$, since $x^2y\notin C$, we get (i), (ii) or (iv). If (i) occurs, then $xyx=xxy$, thus, $xy=yx$, a contradiction. If (ii) or (iv) occurs, then $x(xy)=(xy)^*x=(yx)x$, thus, $(x^2,y)=1$.\cqd

\begin{Lema}Suppose that $\MS$ is anticommutative. Then, for all $y\in G_*$, $x\notin G_*$, and $\alpha$ such that $\alpha y\in \MS$, it holds that:
\begin{itemize}
 \item [(i)] $x^y\in\che x,x^*\chd$.
 \item [(ii)] $xy\in G_*\sse xy\neq yx$.
 \item [(iii)] If $xy\neq yx$, then $\alpha(1+\sigma(x))=0$. If $xy=yx$, then $2\alpha=0$.
 \item [(iv)] $(x,y^2)=(x^2,y)=(xx^*,y)=1$.
 \end{itemize}
\label{as.um.em.cada}
\end{Lema}

\dem In order to prove (i), let a nonzero $\alpha\in R$ such that $\alpha y\in\MS$, then
$$\begin{array}{rcl}
0&=&\alpha y(x+\sigx x^*)+(x+\sigx x^*)\alpha y\\
&=&\alpha yx+\alpha \sigx yx^*+\alpha xy+\alpha\sigx x^*y.
\end{array}$$
Since $\alpha\neq0$ and $x\notin G_*$, $xy\in\che yx, yx^*\chd$, which is equivalent to $x^y\in\che x,x^*\chd$.

%Suppose $xy\neq yx^*$. 

%If $\sigy\neq-1$, then $(xy+\sigxy yx^*)\neq0$ and anticommutes with $(x\inv+\sigma(x\inv)x^{-*})$, so
%$$\begin{array}{rcl}
%0&=&(xy+\sigxy yx^*)(x\inv+\sigma(x\inv)x^{-*})+(x\inv+\sigma(x\inv)x^{-*})(xy+\sigxy yx^*)\\
%&=& xyx\inv+\sigma(x\inv) xyx^{-*}+\sigxy yx^*x\inv+\sigy y+\\
%&&+y+\sigxy x\inv yx^*+\sigma(x\inv)x^{-*}xy+\sigy x^{-*}yx^*.
%\end{array}$$
%Since $\sigy\neq-1$, then $y\in\che xyx\inv,xyx^{-*},yx^*x\inv,x\inv yx^*,x^{-*}xy,x^{-*}yx^*\chd$; on the other hand, since $x\notin G$ and $xy\neq yx^*$, then $y\notin\che xyx^{-*},yx^*x^{-*},x\inv yx^*,x^{-*}xy\chd$, so $y\in \che xyx\inv, x^{-*}yx^*\chd$, which implies $xy=yx$.

%If $\sigy=-1$, then $2y\in\MS$ and it anticommutes with $(x+\sigx x^*)\in\MS$, so
%$$2yx+2\sigx yx^*+2xy+2\sigx x^*y=0,$$
%which implies $xy=yx$ ou $xy=yx^*$. Thus (i) holds.

To prove (ii), it is enough to note that $xy\neq yx\sse xy=yx^*\sse xy\in G_*$.

Let be $\alpha\in R$ such that $\alpha y\in \MS$, thus
$$\begin{array}{rcl}
0&=&\alpha y(x+\sigx x^*)+(x+\sigx x^*)\alpha y\\
&=&\alpha yx+\alpha \sigx yx^*+\alpha xy+\alpha\sigx x^*y.
\end{array}$$
If $xy=yx^*$, in order to guarantee the equality, it is necessary $\alpha(1+\sigx)=0$. In the same way, if $xy=yx$, it is necessary $2\alpha=0$. It proves (iii).

Finally, if $(x,y)\neq1$, then, by item (i),
$$x^2y=xyx^*=y(x^*)^2=yx^2,$$
$$xy^2=yx^*y=y^2x\ \text{and}$$
$$xx^*y=xyx=yx^*x=yxx^*.$$
So (iv) follows.
\cqd

\begin{Lema}Suppose that $\MS$ is anticommutative. Then, for all $x,y\in G_*$ and $\alpha,\beta\in R$ such that $\alpha x,\beta y\in \MS$, it holds
\begin{itemize}
 \item [(i)] $xy=yx\sse xy\in G_*$.
 \item [(ii)] If $xy\neq yx$, then $\alpha\beta=0$; if $xy=yx$, then $2\alpha\beta=0$.
 \item [(iii)] $(x,y^2)=(x^2,y)=1$.
 \end{itemize}
 \label{as.2.simet}
\end{Lema}

\dem Notice that $xy\in G_*$ if, and only if, $xy=(xy)^*=y^*x^*=yx$, since $x,y\in G_*$, thus (i) holds.

If $\alpha,\beta\in R$ satisfies $\alpha x,\beta y\in\MS$, then $\alpha\beta xy=\beta\alpha yx$, then if $xy\neq yx$, we must have $\alpha\beta=0$; if $xy=yx$, then $2\alpha\beta=0$, thus (ii) holds.

Suppose that $xy\neq yx$, so, according to (i), $xy\notin G_*$. Observe that, if $x\in G_*$, then there exists a non zero $\alpha\in R$ such that $\alpha x\in\MS$. In fact, due the compatibility of $*$ and $\sigma$, $x^2=xx^*\in N$, thus, taking $\alpha=(1+\sigx)$, we conclude that $\alpha(1-\sigx)=1-\sigx^2=0$, hence $\alpha x\in \MS$; note that $\alpha\neq 0$ in case that $\sigx\neq-1$, on the other hand, if $\sigx=-1$, we can take $\alpha=2$. So, let $\alpha\in R$ as above, thus, as $\MS$ is anticommutative and $(xy)^*=yx$,
$$\begin{array}{rcl}
0&=&\alpha x(xy+\sigxy yx)+(xy+\sigxy yx)\alpha x\\
&=&\alpha x^2y+\alpha\sigxy xyx+\alpha xyx+\alpha\sigxy yx^2,
\end{array}$$
which implies $x^2y=yx^2$. The other case is analogous.
\cqd

\begin{Lema}If $\MS$ is anticommutative, then $c_x=x^*x\inv\in\centro(G),\ \forall x\in G$.\label{as.cx.central}\end{Lema}

\dem If $x\in G_*$, the result holds trivialy. Let it be $x\notin G_*$ and $y\in G$ such that $(x,y)\neq1$.

If $xy=yx^*$ and $x^*y=yx$, then $x^*x\inv y=x^*yx^{-*}=yxx^{-*}=yx^*x\inv$, since $x^2=(x^*)^2$. So, if $y\notin G_*$ and $xy=yx^*$ or $y\in G_*$, by item (iii) of Lemma \ref{as.2.n.simet} or item (i) of Lemma \ref{as.um.em.cada}, we find $(x^*x\inv,y)=1$.

Suppose that $y\notin G_*$ and $xy\neq yx^*$. By item (i) of Lemma \ref{as.2.n.simet}, $xy\in\che y^*x,x^*y^*\chd$. If $xy=y^*x$, then $y^*x^*=x^*y$ and $x\inv y^*=yx\inv$, thus, $x\inv x^*y=x\inv y^*x^*=yx\inv x^*$. If $xy=x^*y^*$, then $yx=y^*x^*$, thus, $x^*x\inv y^*=x^*x^{-*}y=y=yx^{-*}x^*=y^*x\inv x^*$; in other words, $(x^*x\inv,y^*)=1$, which implies $(x^*x\inv,y)=1$. Therefore $x^*x\inv\in\centro(G)$.
\cqd

The following Lemma is similar to Lemma 3.4 of \cite{JM06}.

\begin{Lema}If $\MS$ is anticommutative, then $x^*=c_xx,\ \forall x\in G$, where $c_ x\in G_*\cap \centro(G)$ and $c_x^2=1$. Furthermore, $c_{xy}=c_xc_y(x,y)$ and, if $(x,y)\neq1$, then $c_{xy}\in\che c_x,c_y,(x,y)\chd$.
\label{as.cxcy}
\end{Lema}

\dem We will prove firstly that $(x,y)\in G_*,\ \forall x,y\in G$. Take $x,y\in G$ such that $(x,y)\neq1$. If $x\notin G_*$ and $y\in G_*$, then, by item (i) of Lemma \ref{as.um.em.cada},  $(x,y)=x\inv x^y=x\inv x^*\in G_*$, since $x^2\in G_*$; the case $x\in G_*$ and $y\notin G_*$ is analogous. If $x,y\notin G_*$ or $x,y\in G_*$, then, by item (ii) of Lemma \ref{as.2.n.simet} or item (i) of Lemma \ref{as.2.simet}, we get $xy, x\inv y\inv\notin G_*$, thus, by item (ii) of Lemma \ref{as.2.n.simet}, $(x,y)=(x\inv y\inv) (xy)\in G_*$ if, and only if, $(xy,x\inv y\inv)=1$; on the other hand, as $x^2, y^2\in \centro(G)$, $(xy,x\inv y\inv)=y\inv x\inv yxxyx\inv y\inv=1$. Thus $(x,y)\in G_*$.

Now we will prove that $(x,y)^2=1,\ \forall x,y\in G$, in other words, $xyx\inv y\inv=yxy\inv x\inv$.

Take $x,y\in G$. If $(x,y)=1$, we have nothing to prove. Suppose that $(x,y)\neq1$. By Lemma \ref{as.2.n.simet} and \ref{as.um.em.cada}, we know that $xy\in\che yx^*,y^*x,x^*y^*\chd$.

If $xy=yx^*$, then, by item (iii) and (i) of Lemma \ref{as.2.n.simet} and \ref{as.um.em.cada}, respectively, we find $yx=x^*y$, so 
$$\begin{array}{rcll}
 xyx\inv y\inv&=&xx^{-*}yy\inv\\
 &=&xx^{-*}\\
 &=&x^*x\inv&(x^2=(x^*)^2)\\
 &=&x^*yy\inv x\inv\\
 &=&yxy\inv x\inv.
 \end{array}$$
The case $xy=y^*x$ is analogous. 

If $xy=x^*y^*$, since $(x,y)\in G_*$ and $x^2,y^2\in \centro(G)$, then 
$$\begin{array}{rcl} 
   xyx\inv y\inv&=&y^{-*}x^{-*}y^*x^*\\
   &=&y\inv x\inv yx\\
   &=&(yy\inv)y\inv (xx\inv)x\inv yx\\
   &=&yy^{-2}xx^{-2}yx\\
   &=&yxy\inv x\inv,
   \end{array}$$
hence $(x,y)^2=1,\ \forall x,y\in G$.

Given $x\in G$, let $c_x=x^*x\inv$ and notice that $c_x\inv=xx^{-*}=x^*x\inv=c_x$, since $x^2\in G_*$, so $c_x^2=1$ and $c_x\in G_*$. Take $y\in G$, then $(xy)^*=c_{xy}xy$ and, on the other hand, by Lemma \ref{as.cx.central}, $(xy)^*=y^*x^*=c_yyc_xx=c_xc_yyx$, in other words, $c_{xy}=c_xc_y(y\inv,x\inv)=c_xc_y(x,y)$.

Notice that Lemmas \ref{as.2.n.simet} and \ref{as.um.em.cada} claim if $(x,y)\neq1$, then one of following holds:
\begin{itemize}
 \item [(a)] $x^*=y\inv xy=yxy\inv$;
 \item [(b)] $y^*=x\inv yx=xyx\inv$;
 \item [(c)] $xy=x^*y^*$ and $xy^*=x^*y$.
\end{itemize}
So if (a) occurs, then $c_x=x\inv x^*=x\inv y\inv xy=(x,y)$ and $c_{xy}=(x,y)c_y(x,y)=c_y$; analogously, if (b) occurs, then $c_y=(x,y)$ and $c_x=c_{xy}$; finally, if (c) holds, then $c_x=x^*x\inv=y^*y\inv=c_y$ and $c_{xy}=(x,y)$, thus (iii) holds.   
\cqd

\begin{Lema}Suppose that $\MS$ is anticommutative. If $x,y\notin G_*$, then:
 \begin{itemize}
\item [(i)] If $(x,y)=1$, then $c_x=c_y$ and $2(1+\sigma(xy))=0=2(\sigx+\sigy)$;
\item [(ii)] If $(x,y)\neq 1$, then: 
 \begin{itemize}
  \item [(a)] $(x,y)=c_x=c_y=c_{xy}$ and $(1+\sigx+\sigy+\sigxy)=0$, or;
  \item [(b)] $(x,y)=c_x\neq c_y=c_{xy}$ and $\sigx=-1$, or;
  \item [(c)] $(x,y)=c_y\neq c_x=c_{xy}$ and $\sigy=-1$, or;
  \item [(d)] $(x,y)=c_{xy}\neq c_x=c_y$ and $\sigxy=-1$ and $\sigx=-\sigy$.
 \end{itemize}
\end{itemize}
\label{as.boa.invol}
\end{Lema}

\dem To show (i) notice that if $(x,y)=1$ and $x,y\notin G_*$, then, by item (iii) of Lemma \ref{as.cxcy}, $c_{xy}=c_xc_y$. By item (ii) of Lemma \ref{as.2.n.simet}, $xy\in G_*$, thus, $c_xc_y=1$, so $c_x=c_y$, since $c_x^2=c_y^2=1$.

To verify (ii), let us consider some cases. If $x,y\in C$, by Theorem \ref{teo.gp14}, and equation (\ref{as.eq.prin}), we obtain (a) easily. If $y\notin C$, then Lemma \ref{as.base.2} implies (a), (b) or (d). If $x\notin C$, using Lemma \ref{as.base.2}, we get (a), (c) or (d).\cqd

The following proposition will be usefull in the classification of the group rings such that $\MS$ is anticommutative, as well it gives a better description of the involution of item (B) of Theorem \ref{teo.gp14}.

\begin{Prop} If $\MS$ is anticommutative, then one of the following holds:
\begin{itemize}
 \item [(i)] $G$ is abelian and $*=Id$.
 \item [(ii)] $G$ is abelian, $*\neq Id$, and exists $s\in G_*\cap\centro(G)$ such that $x^*=\che x,xs\chd,\ \forall x\in G$.
 \item [(iii)] $G'=\che 1,s\chd\subset G_*\cap\centro(G)$, $s^2=1$ and $x^*=\che x,xs\chd,\ \forall x\in G$.
   \item [(iv)] $*|_C\neq Id$, $*|_N=Id$, $G'=\lan s,t\ran\subset G_*\cap\centro(G)$, where $s$ and $t$ are the only nontrivial commutators of $G$, $C'=\che1\chd$ or $C'=\che 1,s\chd$, $s^2=t^2=1$, $x^*=\che x,xs\chd \forall x\in C$, and $x^*=\che x,xt\chd\ \forall x\notin C$. Futhermore, $(G\menos (C\cup G_*),G\menos (C\cup G_*))\subset\che 1,t\chd$, $(C\menos G_*,G)\subset\che 1,s\chd$, and $(C\menos G_*,G\menos (C\cup G_*))=\che s\chd$.
 \end{itemize}
\label{as.desc.g}
 \end{Prop}

\dem Suppose that $G$ is abelian and $*\neq Id$. Fix $x\notin G_*$ and denotes $c_x=s$. Given $y\notin G_*$, as $G$ is abelian, by item (i) of Lemma \ref{as.boa.invol}, $c_x=c_y$, thus $c_y=s$ and $y^*=ys,\ \forall y\notin G_*$, hence (ii) holds.

Applying Theorem \ref{teo.gp14} and what we prove above, we can conclude that the involution in $C$ is given by $x^*\in\che x,sx\chd,\ \forall x\in C$, for some $s\in G_*$.

We can assume that $G$ is non abelian.

Suppose that $*|_C=Id$, thus $C$ is abelian. Fix $z\notin G_*$ and $s=c_z$, thus, $z\notin C$. Let $x\notin G_*$, and applying Lemma \ref{as.boa.invol}, we have $(z,x)=1$ and $c_x=s$, or, $(z,x)=s=c_x=c_{zx}$, otherwise, $*|_C\neq Id$; thus, $x^*\in\che x,sx\chd\ \forall x\in G$, and $(x,y)\in\che 1,s\chd,\ \forall x,y\notin G_*$. Let $x,y\in G$ such that $(x,y)\neq1$. If $x\notin G_*$ and $y\in G_*$, by Lemma \ref{as.um.em.cada}, we have $xy\in G_*$, then, by item (ii) of Lemma \ref{as.cxcy}, $c_{xy}=c_xc_y(x,y)$; in other words, $1=c_x(x,y)$, which implies $(x,y)=s$, for $c_x^2=1$. If $x,y\in G_*$, by Lemma \ref{as.2.simet}, $xy\notin G_*$, so $c_{xy}=s$, then, by item (ii) of Lemma \ref{as.cxcy}, $c_{xy}=c_xc_y(x,y)$, then, $s=(x,y)$. This way, we conclude that (iii) holds.

Suppose $*|_C\neq Id$. Let $x\in C\menos G_*$ and $s=c_x$ given by Theorem \ref{teo.gp14}. let $y\in G$ such that $(x,y)\neq1$. Suppose $y\notin G_*$; if $y\in C$, by Theorem \ref{teo.gp14}, $(x,y)=s$; if $y\notin C$, by Lemma \ref{as.boa.invol}, (a) or (b) holds, since $x\in C$ and $y\notin C$, thus, in both cases, $(x,y)=c_x=s$. On the other hand, if $y\in G_*$, by Lemma \ref{as.um.em.cada}, $xy\in G_*$, then, by item (ii) of Lemma \ref{as.cxcy}, $c_{xy}=c_xc_y(x,y)$, so, $1=c_x(x,y)$, which implies $(x,y)=s$, since $s^2=1$. Thus $(C\menos G_*,G)\subset\che 1,s\chd$, since $c_y=s\ \forall y\in C\menos G_*$.

%Vamos mostrar que existe $x\notin G_*$ tal que $x\notin\centro(G)$.

%Suponha por absurdo que $G\menos G_*\subset\centro(G)$. Seja $x\notin G_*$ e note que $xy\notin G_*,\ \forall y\in G_*$. Dados $y,z\in G_*$, temos que $(z,y)=(xz,y)=1$, já que $x,xz\in\centro(G)$, portanto, $G$ é abeliano, contradição.

%Sejam $x\notin G_*$ tal que $x\notin\centro(G)$ e $y\in G$ tal que $(x,y)\neq1$. Se $y\notin G_*$, então, pelo item (ii) do Lema \ref{as.boa.invol}, temos que $c_x,c_y$ ou $c_{xy}=s$; se $y\in G_*$, então pelo Lema \ref{as.um.em.cada}, temos que $x^*=x^y=xs$, ou seja, $c_x=s$. Portanto podemos afirmar que existe $x\notin G_*$ tal que $c_x=s$.

Let $x,y\in G\menos (G_*\cup C)$ and $c_x=t$. By Lemma \ref{as.boa.invol}, if $(x,y)=1$, then $t=c_x=c_y$; if $(x,y)\neq 1$, then $t=c_x=c_y=(x,y)$, or, $t=c_x=c_y$, $\sigma(xy)=-1$, $s=c_{xy}=(x,y)$. Thus, the involution in $G\menos C$ is given by $x^*=\che x,xt=tx\chd\ \forall x\notin C$. Let $x\in G\menos (G_*\cup C)$ and $y\in G_*$ such that $(x,y)\neq1$; this way, by Lemma \ref{as.um.em.cada}, $xy\in G_*$, so by Lemma \ref{as.cxcy}, $c_{xy}=c_xc_y(x,y)$, in other words, $1=c_x(x,y)$, then $(x,y)=t$, thus $(G\menos G_*,G)\subset\che 1,s,t\chd$. Finnaly, let $x,y\in G_*$, such that $(x,y)\neq1$; by Lemma \ref{as.2.simet}, $xy\notin G_*$ and again by Lemma \ref{as.cxcy}, $c_{xy}=c_xc_y(x,y)$, so, $(x,y)=c_{xy}=s$ or $t$; then we conclude that $G$ has at most two nontrivial commutators $s$ e $t$; besides that $x^*=\che x,xs\chd,\ \forall x\in C$ and $x^*=\che x,xt\chd\ \forall x\notin C$. 

If $t=s$, (iii) holds. Suppose $t\neq s$, thus, $G\menos (G_*\cup C)\neq\emptyset$. If $*|_N\neq Id$, then taking $y\in N\menos G_*$ and $x\notin (G_*\cup C)$, by Lemma \ref{as.boa.invol}, we find $s=c_x=c_y=t$, a contradiction, since $s\neq t$, thus $*|_N= Id$. Finnaly, if exist $x\in C\menos G_*$ and $y\in G\menos(C\cup G_*)$ such that $(x,y)=1$, again, by Lemma \ref{as.boa.invol}, $s=c_x=c_y=t$, a contradiction; moreover, as $(G\menos G_*,G)\subset\che 1,s\chd$, so we conclude that $(C\menos G_*,G\menos(C\cup G_*))=\che s\chd$.
\cqd

%Por outro lado se dado $x\in C\menos G_*$, temos pelo Lema \ref{as.boa.invol} que $c_x=t$ ou $s$. Como mostramos que existe algum elemento em $G$ tal que $c_x=s$, se $t=s$, temos que (a) ocorre, caso $t\neq s$, então vale (b).

%{\bf Afirmação:} Existe $x\notin G_*$ tal que $x\notin\centro(G)$.

%Suponha por absurdo que $G\menos G_*\subset\centro(G)$. Seja $x\notin G_*$ e note que $xy\notin G_*,\ \forall y\in G_*$. Dados $y,z\in G_*$, temos que $(z,y)=(xz,xy)=1$, já que $xy,xz\in\centro(G)$, portanto, $G$ é abeliano, contradição.

%Sejam $x\notin G_*$ tal que $x\notin\centro(G)$ e $y\in G$ tal que $(x,y)=s$. Se $y\notin G_*$, então, pelo item (ii) do Lema \ref{as.boa.invol}, temos que $c_x,c_y$ ou $c_{xy}=s$; se $y\in G_*$, então pelo Lema \ref{as.um.em.cada}, temos que $x^*=x^y=xs$, ou seja, $c_x=s$. Portanto podemos afirmar que existe $x\notin G_*$ tal que $c_x=s$.

%Seja $y\in G\menos (G_*\cup C)$ e $c_y=t$, temos que, pelo Lema \ref{as.boa.invol}, se $z\in G\menos (G_*\cup C)$, se $(y,z)=1$, então $t=c_y=c_z$ e caso $(y,z)=s$, então $t=c_y=c_z=s$ ou $t=c_y=c_z$ e $\sigma(yz)=-1$. Portanto a involução em $G\menos C$ é dada por $x^*=\che x,xt=tx\chd\ \forall x\notin C$. Por outro lado se dado $x\in C\menos G_*$, temos pelo Lema \ref{as.boa.invol} que $c_x=t$ ou $s$. Como mostramos que existe algum elemento em $G$ tal que $c_x=s$, se $t=s$, temos que (a) ocorre, caso $t\neq s$, então vale (b).

\pagebreak

\begin{Teo} Let $R$ be a commutative ring of $char(R)\neq 2$ and $G$ a group with an involution $*$ compatible with an orientation $\sigma$. $\MS$ is anticommutative if, and only if, the following holds:
\begin{itemize}
 \item [(i)] 
 \begin{itemize}
  \item [(a)] $char(R)=8,4$, $G$ is abelian and $*=Id$; or
  \item [(b)] $char(R)=4$, $*\neq Id$ and $G$ satisfies one of the following: 
  \begin{itemize}
   \item [(1)] $G$ is abelian and exists $s\in N_*\cap\centro(G)$ such that $x^*=\che x,xs\chd,\ \forall x\in G$.
   \item [(2)] $G'=\che 1,s\chd\subset N_*\cap\centro(G)$, $s^2=1$, and  $x^*=\che x,xs\chd,\ \forall x\in G$.
   \item [(3)] $*|_C\neq Id$, $*|_N=Id$, $G'=\lan s,t\ran\subset N_*\cap\centro(G)$, where $s$ and $t$ are the only nontrivial commutators of $G$, $C'=\che1\chd$ or $C'=\che 1,s\chd$, $s^2=t^2=1$, $x^*=\che x,xs\chd,\ \forall x\in C$, and $x^*=\che x,xt\chd\ \forall x\notin C$. Moreover, $(G\menos (C\cup G_*),G\menos (C\cup G_*))\subset\che 1,t\chd$, $(C\menos G_*,G)\subset\che 1,s\chd$ e $(C\menos G_*,G\menos (C\cup G_*))=\che s\chd$.
  \end{itemize}
  \end{itemize}
\item [(ii)] If $x,y\notin G_*$ and $(x,y)=1$, then $2(1+\sigma(xy))=0=2(\sigx+\sigy)$; if $\sigx\neq-1\neq\sigy$, and $(x,y)\neq1$, then $(1+\sigx+\sigy+\sigxy)=0$.
 \item [(iii)] $\forall x\notin G_*,\ y\in G_*$ and $\alpha\in R$ such that $\alpha y\in \MS$, if $xy\neq yx$, then $\alpha(1+\sigma(x))=0$; if $xy=yx$, then $2\alpha=0$.
 \item [(iv)] $\forall x,y\in G_*$ and $\alpha,\beta\in R$ such that $\alpha x,\beta y\in \MS$. If $xy\neq yx$, then $\alpha\beta=0$; if $xy=yx$, then $2\alpha\beta=0$.
 \end{itemize}
 \label{as.teo.1}
\end{Teo}

\dem Suppose that $\MS$ is anticommutative. By Lemma \ref{as.char4.char8}, we know that $char(R)=4$ or $8$. We can easily check that $Id$ defines an involution only if $G$ is abelian, thus (i)(a) occurs.

Suppose that $*\neq Id$. Thus, the Lemma \ref{as.*.neq.id} implies $char(R)=4$ and, Proposition \ref{as.desc.g} ensure (i)(b), so (i) holds. 

The item (ii) follows by Lemma \ref{as.cxcy}, item (iii) by item (iii) of Lemma \ref{as.um.em.cada}, and item (iv) by item (ii) of Lemma \ref{as.2.simet}.

To prove the converse, denotes $c_a=a^*a\inv,\ \forall a\in G$ and notice that, by (i), $c_a\in\centro(G)$.

Let $x,y\notin G_*$; we must prove that 
 $$(x+\sigx x^*)(y+\sigy y^*)+(y+\sigy y^*)(x+\sigy x^*)=0,$$ in other words,
\begin{equation}xy+yx+\sigma(xy)xyc_xc_y+\sigma(xy)yxc_xc_y+\sigma(y)xyc_y+\sigma(y)yxc_y+\sigma(x)yxc_x+\sigma(x)xyc_x,
\label{as.eq.final}
\end{equation}
vanishes.

Notice that (a), (1), (2), and condition $(C\menos G_*,G\menos (C\cup G_*))$ in (3), ensure that if $(x,y)=1$, then $c_x=c_y$ and, by (ii), $2(1+\sigxy)=0=2(\sigx+\sigy)$, thus (\ref{as.eq.final}) could be rewriten by
$$xy+xy+\sigma(xy)xy+\sigma(xy)xy+\sigma(y)xyc_x+\sigma(y)xyc_x+\sigma(x)xyc_x+\sigma(x)xyc_x,$$
in other words,
$$2(1+\sigxy)xy+2(\sigx+\sigy)xyc_x=0.$$

Suppose that $(x,y)\neq1$. If $\sigx=-1$, then, by (i), $c_x=s=(x,y)$, thus (\ref{as.eq.final}), becomes
$$\begin{array}{c}
   xy+xys+\sigma(xy)xysc_y+\sigma(xy)xys^2c_y+\sigma(y)xyc_y+\sigma(y)xysc_y+\sigma(x)xys^2+\sigma(x)xys\\
   =\\
   (1+\sigx)(xy+xys)+(\sigy+\sigxy)(xyc_y+xysc_y)\\
   =\\
   0;
\end{array}$$
if $\sigy=-1$, we can prove it in an analogous way; if $\sigx\neq-1\neq\sigy$, then, if occurs, we obtain $c_x=c_y=s=(x,y)$; if (3) occurs, then $x,y\notin C$, since $*|_N=Id$, thus $c_x=c_y=t=(x,y)$. In both cases,
$$\begin{array}{c}
   xy+xyc_x+\sigma(xy)xyc_x^2+\sigma(xy)xyc_x^3+\sigma(y)xyc_x+\sigma(y)xyc_x^2+\sigma(x)xyc_x^2+\sigma(x)xyc_x\\
   =\\
   (1+\sigx+\sigy+\sigxy)(xy+xyc_x)\\
   =\\
   0.
\end{array}$$

Let $x\notin G_*$ and $y\in G_*$. we must prove that $(x+\sigx x^*)$ and $\alpha y$ anticommutes, where $\alpha\in R$ is such that $\alpha y\in \MS$. If $(x,y)\neq 1$, then, by Lemma \ref{as.um.em.cada}, $xy=yx^*$, $x^*y=yx$, and, by (iii), $\alpha(1+\sigx)=0$, so
$$\begin{array}{rcl}
(x+\sigx x^*)\alpha y+\alpha y(x+\sigx x^*)&=&\alpha xy+\alpha \sigx x^*y+\alpha yx+\alpha \sigx yx^*\\
&=&\alpha xy+\alpha \sigx yx+\alpha yx+\alpha \sigx xy\\
&=&\alpha(1+\sigx)(xy+yx)\\
&=&0.
\end{array}$$
If $(x,y)=1$, then $(x^*,y)=1$ and we can proceed in an analogous way to prove that $(x+\sigx x^*)$ anticommutes with $\alpha y$, since, $2\alpha=0$, by (iii).

Let $x,y\in G_*$ and $\alpha,\beta\in R$ such that $\alpha x,\beta y\in\MS$. Using item (iv), we can prove as in the previous case, that $\alpha x$ anticommutes with $\beta y$.
\cqd

%\section{Skew-symmetric Elements Anticommute}
%\input{2.anticomutatividade.antissimetricos.tex}

%%%%%%%%%% BIBLIOGRAFIA %%%%%%%%%%
%

\renewcommand{\bibname}{Referências}
\addcontentsline{toc}{chapter}{\bibname}       % Adiciona ao sumário a referência à Bibliografia.

%\addcontentsline{toc}{chapter}{\indexname}    % Adiciona ao sumário a referência ao Índice Remissivo.

%\printindex
\end{document}